\documentclass[12pt,a4paper]{article}
\usepackage{url}
\usepackage{amsmath}    
\usepackage{amssymb}
\usepackage{graphicx}   
\usepackage{verbatim}   
\usepackage{color}      
\usepackage{hyperref}   
\usepackage{enumerate}
\usepackage{listings}
\usepackage{xcolor}
\usepackage{slashed}
\usepackage{amsmath}    
\usepackage{amssymb}
\usepackage{mathrsfs}
\usepackage{bbold}
\usepackage[makeroom]{cancel}

\usepackage[final]{showlabels}

\usepackage{amsthm}
\newtheorem{remark}{Remark}

\renewcommand{\title}{Truly multi-dimensional all-speed methods for the Euler equations}

\newcommand{\authorOne}{Wasilij Barsukow\footnote{Bordeaux Institute of Mathematics, Bordeaux University and CNRS/UMR5251, Talence, 33405 France}}

\begin{document}

\newcommand{\del}{\partial}
\renewcommand{\theta}{\vartheta}
\renewcommand{\phi}{\varphi}
\newcommand{\vecc}[2]{\left ( \begin{array}{c}#1\\#2\\ \end{array}\right )}
\newcommand{\veccc}[3]{\left ( \begin{array}{c}#1\\#2\\#3\\ \end{array}\right )}
\newcommand{\dd}{\mathrm{d}}
\newcommand{\ee}{\mathrm{e}}
\newcommand{\ii}{\mathrm{i}}
\newcommand{\id}{\mathbb{1}}
\newcommand{\atanh}{\,\text{artanh}\,}
\newcommand{\atan}{\arctan}
\newcommand{\back}{\!\!\!}
\newcommand{\nicefrac}[2]{#1 / #2}
\newcommand{\bboxed}[1]{\text{\textsc{Conjecture}: }\boxed{\boxed{#1}}}
\renewcommand{\and}{\wedge}
\newcommand{\primitive}{\text{\textsc{primitive}}}
\newcommand{\lint}{\int\limits}
\renewcommand{\div}{\mathrm{div\,}}
\renewcommand{\vec}{\mathbf}
\newcommand{\todo}[1]{{\color{red}#1}}
\newcommand{\new}[1]{{\color{blue}#1}}

\begin{center} \Large
\title

\vspace{1cm}

\date{}
\normalsize

\authorOne
\end{center}

\begin{abstract}
Several recent all-speed time-explicit numerical methods for the Euler equations on Cartesian grids are presented and their properties assessed experimentally on a complex application. These methods are truly multi-dimensional, i.e. the flux through an interface also depends on the values in cells adjacent to the endpoints of the edges (corners). 

\end{abstract}
\section{Introduction}

The rescaled Euler equations read
\begin{align}
\del_t \rho + \nabla \cdot (\rho \vec v) &= 0   & \rho & \colon \mathbb R^+_0 \times \mathbb R^d \to \mathbb R^+, \label{eq:eulercons1}\\
\del_t (\rho \vec v) + \nabla \cdot (\rho \vec v \otimes \vec v + \frac{p}{\epsilon^2} \id) &= 0 & \vec v & \colon \mathbb R^+_0 \times \mathbb R^d \to \mathbb R^d, \\
\del_t e + \nabla \cdot ((e + p) \vec v) &= 0 & e,p & \colon \mathbb R^+_0 \times \mathbb R^d \to \mathbb R^+, \label{eq:eulercons3} \\
\frac{p}{\gamma-1} + \frac12 \epsilon^2 \rho |\vec v|^2 &= e ,
\end{align}
where $\epsilon \to 0$ corresponds to the limit of low Mach number. This is a singular limit: for well-prepared initial data, the solutions tend to those of the incompressible Euler equations as $\epsilon \to 0$. 

Numerical diffusion that is used to stabilize time-explicit discretizations of \eqref{eq:eulercons1}--\eqref{eq:eulercons3} is often introduced by means of a Riemann solver, and thus inspired by supersonic phenomena. In subsonic and low Mach number flow, this diffusion is excessive and dominates the numerical solution unless the grid is refined strongly (see e.g. \cite{guillard04}). The excessive diffusion is also worrying because it points to a fundamental deficiency of Godunov-type methods (see e.g. \cite{barsukow17}). 

The suggestions to overcome this problem traditionally amount to 
\begin{enumerate}[1.]
\item removing certain terms from the diffusion matrix (or to multiply them by a factor $f \in \mathcal O(\epsilon)$) (\textbf{low Mach fixes}, e.g. \cite{dellacherie10}), or to 
\item using central differences for at least the acoustic operator, which necessitates (partially) implicit time integration (\textbf{IMEX}, e.g. \cite{boscheri21a}).
\end{enumerate}

Time explicit methods are easy to implement and are computationally inexpensive, but the CFL condition requires $\Delta t \in \mathcal O(\epsilon)$. Provided that the implicit methods are well-conditioned, as $\epsilon \to 0$, they will always outperform explicit methods in terms of computational efficiency. For good and well-implemented methods, the break-even point is typically at $\epsilon \sim 10^{-2}$.

This paper focuses on time-explicit methods for the following reason: If very low values of $\epsilon$ are to be used, and acoustic phenomena (sound waves, density stratifications, etc.) are of no interest, then it might be more practical to use an incompressible code. If the acoustic phenomena need to be resolved, out of considerations of accuracy, even for an implicit method the time step would be chosen close to that of an explicit method (in view of the above-mentioned break-even point, even choosing the time step to be 100 times that of an explicit method would barely justify using an implicit method).

Two strategies that reduce the numerical diffusion while retaining enough of it to allow explicit integration in time, and avoiding the introduction of ad-hoc factors have been developed recently:
\begin{enumerate}[1.] \setcounter{enumi}{2}
\item In the methods from \cite{barsukow20cgk}, the divergence $\nabla \cdot \vec v \in \mathcal O(\epsilon)$ plays the role of the factor $f \in \mathcal O(\epsilon)$, because it seems to arise naturally in just the right place.
\item \cite{barsukow21yee} points out that, in fact, central differences can be stable upon \emph{explicit} time integration for certain systems, e.g. the acoustic equations.
\end{enumerate}

After some preliminary remarks concerning notation (Section \ref{sec:notation}) as well as the advective operator (Section \ref{sec:advection}) these two strategies shall be presented (Sections \ref{sec:cgk} and \ref{sec:yee}) and compared to each other (Section \ref{sec:numerical}). The numerical methods are considered on two-dimensional Cartesian grids, with cells $[x_{i-\frac12}, x_{i+\frac12}] \times [y_{j-\frac12}, y_{j+\frac12}]$.

\section{Notation for finite difference operators} \label{sec:notation}

In order to shorten finite difference formulae, and also in order to make them better readable, the following short-hand notation is used:
\begin{align}
 [a]_{i+\frac12} := a_{i+1} - a_i \quad \text{the jump/difference operator,} \\
 \{a\}_{i+\frac12} := a_{i+1} + a_i \quad \text{the summation operator.}
\end{align}
Both are centered at a cell interface, but one also has the cell-centered versions
\begin{align}
 [b]_{i\pm\frac12} &:= b_{i+\frac12} - b_{i-\frac12}, &
 [a]_{i\pm1} &:= a_{i+1} - a_{i-1}
\end{align}
and similarly for the sums. This notation can be combined in one dimension, e.g. $[[a]]_{i\pm\frac12} := [a]_{i+\frac12} - [a]_{i-\frac12} = a_{i+1} - 2 a_i + a_{i-1}$ and $\{ [a] \}_{i\pm\frac12} = [a]_{i\pm1}$
as well as in multiple dimensions (making sure not to interchange the order of the indices)
\begin{align}
 [a_i]_{j+\frac12} \equiv [a]_{i,j+\frac12} &:= a_{i,j+1} - a_{ij}, \\
 [a]_{i+\frac12,j} \equiv [a_{\cdot, j}]_{i+\frac12} &:= a_{i+1,j} - a_{ij}, \\
 \{[a]_{i+\frac12} \}_{j+\frac12} &= [a]_{i+\frac12,j+1} + [a]_{i+\frac12,j} \quad \text{etc.}
\end{align}

\begin{remark}
 The number of square brackets gives the order of the derivative operator that is approximated by the given finite difference formula. The curly brackets do not contain any hidden factors, an average therefore is $\frac12 \{ \cdot \}$.
\end{remark}

\section{Upwinding of the advection operator and Lagrange-Projection methods} \label{sec:advection}

The Euler equations can be considered as made up of three operators in the following way (denote in 2-d $\vec v = (u, v)$):
\begin{align}
 \del_t \rho &&&+ (u \del_x + v \del_y) \rho &&+ \rho (\del_x u + \del_y v) &&= 0 \label{eq:euler1}\\
 \del_t (\rho u) &&&+ (u \del_x + v \del_y) (\rho u) &&+ \rho u (\del_x u + \del_y v) &&+ \frac{\del_x p}{\epsilon^2} = 0\\
 \del_t (\rho v) &&&+ (u \del_x + v \del_y) (\rho v) &&+ \rho v (\del_x u + \del_y v) &&+ \frac{\del_y p}{\epsilon^2} = 0\\
 \del_t e &&&+ (u \del_x + v \del_y) e &&+ e (\del_x u + \del_y v) &&+ \del_x (up) + \del_y (vp)= 0 \label{eq:euler4}\\
 \nonumber &&&\text{\,\,\,\,\,\,advection}&&\text{\,\,\,\,\,\,compression}&&\text{nonlinear acoustics}
\end{align}
where advection $+$ compression is the pressureless Euler system. For stability, numerical diffusion needs to be added: upwinding for advection  and characte\-ristics-based upwinding for acoustics, or possibly more complicated multi-dimensional operators as discussed later. This Section deals with the question whether the compressive terms should be taken into account by central differences, or whether some kind of diffusion associated to these terms needs to be included as well. Essentially, there are two choices that can be found amongst the various numerical methods in the literature. The prototyic equation to study is that of advection
\begin{align}
 \del_t q + \del_x (U q) &= 0 \label{eq:advectionnconst}
\end{align}
with a non-constant speed $U(x)$. The first choice ignores the fact that $U$ varies (in space) and the discretization (for $U>0$)
\begin{align}
 q_i^{n+1} &= q_i^n - \frac{\Delta t}{\Delta x} \left( U(x_{i+\frac12}) q_i^n- U(x_{i-\frac12}) q_{i-1}^n \right ) 
\end{align}
of \eqref{eq:advectionnconst} has the same numerical diffusion as that of
\begin{align}
 \del_t q + U \del_x q &= 0
\end{align}
i.e. they differ only in the central part. This would be the case for the Roe solver, for example, and also for the methods proposed in \cite{leveque02}, Section 9.2.

Lagrange-Projection methods would track the volume contained between two characteristics and would arrive at the following method
\begin{align}
 q_i^{n+1} &= q_i^n - \frac{\Delta t}{\Delta x} \left( U(x_{i+\frac12}) \frac{q_i^n}{1 + \Delta t \frac{U(x_{i+\frac12}) - U(x_{i-\frac12})}{\Delta x}} \right.  \label{eq:lagrprojconsadv} \\\nonumber &\phantom{mmmmmmmmmm} \left. - U(x_{i-\frac12}) \frac{q_{i-1}^n}{1 + \Delta t \frac{U(x_{i-\frac12}) - U(x_{i-\frac32})}{\Delta x}}  \right ) 
\end{align}
easily recognizable by the presence of denominators $\frac{1}{1 + \Delta t \del_x U}$. Relaxation solvers such as \cite{berthon06,bouchut04,bouchut09,chalons10,girardin14} contain similar denominators as well. Further details on these different approaches can be found in \cite{barsukow20cgk,barsukow21yee}. There, it is generally found to be beneficial for stability to include the denominator.

\section{All-speed extension to multiple dimensions} \label{sec:cgk}

The key observation of \cite{barsukow20cgk} is that the diffusive terms that are prohibiting low Mach compliance often are of the form $\del_x(c_1 \del_x u)$ in the equation for the $x$-momentum $\rho u$ and $\del_y (c_2 \del_y v)$ in the equation for the $y$-momentum $\rho v$, with $c_1, c_2$ being terms that depend on the method. While low Mach fixes would endow these two terms with prefactors $f \in \mathcal O(\epsilon)$, the strategy advocated in \cite{barsukow20cgk} is to extend them to $\del_x(c_1 (\del_x u + \del_y v))$ and $\del_y (c_2 (\del_x u + \del_y v)$, respectively. They then contain the divergence, which in the limit $\epsilon \to 0$ would become $\mathcal O(\epsilon)$, thus achieving the same as the prefactor. Due to the appearance of the mixed derivative $\del_x \del_y$ the resulting method necessarily is truly multi-dimensional. The advantage of this approach is that the one-dimensional method (with possibly all its good properties concerning stability, entropy, etc.) remains unmodified, that the method is all-speed, and that the low Mach compliance is achieved without ad-hoc factors.

As now the divergence would appear in both the central part and the diffusion, its respective discretizations need to be related, such that the cell-centered diffusion contains the difference of two edge-centered divergences and such that the central divergence is their average. Finally, due to the fact that there is a diffusion in $x$-direction and one in $y$-direction, the fundamental object is, in fact, a node-based divergence. A divergence centered on an edge is obtained as the average of the two divergences at the nodes. It has been shown in \cite{barsukow17a} that this is the unique procedure for symmetric discretizations.

The strategy of multi-dimensional all-speed extension, as explained above, has been applied to two methods in \cite{barsukow20cgk}. First, the Lagrange-Projection-type method from \cite{chalons13} was considered whose numerical flux in 1-d is
\begin{align}
 \hat f_{i+\frac12} := \left (  0 , \frac{p^*_{i+\frac12}}{\epsilon^2} , 0 , p^*_{i+\frac12} u^*_{i+\frac12}  \right) + u^*_{i+\frac12} \left\{  \frac{Q^n}{L}  \right\}_{i+\frac12}  - |u^*_{i+\frac12}| \left[  \frac{Q^n}{L}  \right]_{i+\frac12} ,
\end{align}
having defined
\begin{align}
 u^*_{i+\frac12} &:= \frac{ \{ u \}_{i+\frac12} }{2} - \frac1{2a \epsilon} [p]_{i+\frac12} & 
 p^*_{i+\frac12} &:= \frac{\{p \}_{i+\frac12} }{2} - \frac{a \epsilon}{2} [u]_{i+\frac12}, \label{eq:1dlagprojrelaxacp}
\end{align}
\vspace{-.4cm}
\begin{align}
 Q_i^n &:= \left(\rho_i^n, (\rho u)^n_i - \frac{\Delta t}{\Delta x} \frac{1}{\epsilon^2} [p^*]_{i\pm\frac12}, (\rho v)^n_i, e_i^n - \frac{\Delta t}{\Delta x} [u^* p^*]_{i\pm\frac12} \right), \\
 L_i &:= 1 + \frac{\Delta t}{\Delta x} [u^*]_{i\pm\frac12}.
\end{align}

\begin{remark}
 The quantities $Q$ are conserved quantities evolved in time through just the acoustic operator. The flux contains upwinding (visible in the term containing $|u^*|$) and also accounts additionally for the compression terms through division by $L$. $a$ is the relaxation speed in the relaxation solver associated with the acoustic sub-system, and fulfills $a > \max (\rho c)$.
\end{remark}

Observe how the diffusion of the momentum equation takes its origin in $p^*$, and more precisely in the term $[u]_{i+\frac12}$. It is this term that the authors of \cite{chalons13} propose to multiply by $\epsilon$ in order to achieve low Mach compliance. Here, instead, expressions \eqref{eq:1dlagprojrelaxacp} shall be replaced by
\begin{align}
 u^*_{i+\frac12,j} &=  \frac{\{\{ \{ u \}_{i+\frac12} \}\}_{j\pm\frac12} }{8} - \frac1{2a\epsilon} \frac{\{\{ [p]_{i+\frac12} \}\}_{j\pm\frac12}}{4}, \label{eq:lagprojmultidfirst}\\
 p^*_{i+\frac12,j} &= \frac{\{\{\{p \}_{i+\frac12} \}\}_{j\pm\frac12} }{8} - \frac{a\epsilon}{2}\left( \frac{\{\{ [u]_{i+\frac12} \}\}_{j\pm\frac12}}{4} + \frac{\Delta x}{\Delta y} \frac{[\{v\}_{i+\frac12}]_{j\pm1} }{4}\right ), \label{eq:lagprojmultidfirstpxdir} \\
 v^*_{i,j+\frac12} &=  \frac{ \{\{ \{ v \}\}_{i\pm\frac12}\}_{j+\frac12} }{8} - \frac1{2a\epsilon} \frac{ [ \{\{ p \}\}_{i\pm\frac12}]_{j+\frac12}}{4}. \label{eq:lagprojmultidfirst2}
\end{align}
It is also natural to redefine
\begin{align}
 L_{ij} &:= 1 + \frac{\Delta t}{\Delta x} [u^*]_{i\pm\frac12,j} + \frac{\Delta t}{\Delta y} [v^*]_{i,j\pm\frac12}, \label{eq:newL}\\
 (Q_{ij}^n)_e &:= e^n_{ij} - \frac{\Delta t}{\Delta x} [u^* p^*]_{i\pm\frac12,j} - \frac{\Delta t}{\Delta y} [v^* p^*]_{i,j\pm\frac12} , \label{eq:lagprojmultidlast}
\end{align}
where $(Q_{ij}^n)_e$ is the energy-component of the vector $Q_{ij}^n$. This method shall be referred to as \textbf{Method A}. Now, an edge-based divergence is appearing in \eqref{eq:lagprojmultidfirstpxdir} (the average of two node-based divergences), while the divergence in \eqref{eq:newL} turns out to be the average of four such edge-based divergences. They thus all vanish simultaneously, which is one of the ingredients of the proof of low Mach compliance. For further details, see \cite{barsukow20cgk}.

For the relaxation solver from \cite{chalons10}, the other method investigated in \cite{barsukow20cgk}, one finds for the intermediate states
\begin{align}
 u^*_{i+\frac12} &= \frac{\{ u \}_{i+\frac12}}{2}  - \frac{1}{2a \epsilon} [p]_{i+\frac12}, & 
 p^*_{i+\frac12} &= \frac{ \{ p \}_{i+\frac12} }{2} - \frac{a \epsilon}{2} [u]_{i+\frac12}, \label{eq:relax1dpstar}\\
 \rho^*_{i+\frac12,\text L} &=  
 \frac{\rho_i}{1 + \rho_i \epsilon \frac{[u]_{i+\frac12}}{2a} - \rho_i\frac{[p]_{i+\frac12}}{2 a^2}}\label{eq:relax1drhoL}, &
 \frac{e^*_{i+\frac12,\text L} }{ \rho^*_{i+\frac12,\text L}} &= \frac{e_i}{\rho_i} + \epsilon \frac{p_i u_i - p^*_{i+\frac12} u^*_{i+\frac12}}{a} \quad \text{etc.}
\end{align}
Equation \eqref{eq:relax1dpstar} is the same as before, up to the definition of $a$ which now is more complicated (see \cite{chalons10} for more details). The flux associated with the intermediate state is (if $u^*_{i+\frac12} > 0$)
\begin{align}
 f_{i+\frac12} = \left (                      u^*_{i+\frac12} \rho^*_{i+\frac12,\text L} ,
                          \rho^*_{i+\frac12,\text L} (u^*_{i+\frac12})^2 + \frac{p^*_{i+\frac12}}{\epsilon^2} ,
                         \rho^*_{i+\frac12,\text L} u^*_{i+\frac12} v_i ,
                          u^*_{i+\frac12} (e^*_{i+\frac12,\text L}  + p^*_{i+\frac12})
                         \right ). \label{eq:relax1dflux}
\end{align}

In \cite{barsukow20cgk}, therefore, equations \eqref{eq:relax1dpstar} are again replaced by \eqref{eq:lagprojmultidfirst}--\eqref{eq:lagprojmultidfirst2}, and the intermediate density is also naturally replaced by
\begin{align}
 \rho^*_{i+\frac12,j,\text{L}} &= 
 \frac{\rho_{ij}}{1 + \frac{\rho_{ij} \epsilon}{2a} \left( \frac{\{\{[u]_{i+\frac12}\}\}_{j\pm\frac12}}{4} + \frac{\Delta x}{\Delta y} \frac{[\{  v\}_{i+\frac12}]_{j\pm1}}{4} \right ) - \frac{\rho_{ij}}{ 2 a^2 } \frac{\{\{ [p]_{i+\frac12} \}\}_{j\pm\frac12}}{4}} \label{eq:relax2drhoL}
\end{align}
and similarly for $\rho^*_{i+\frac12,j,\text{R}}$. This method shall be referred to as \textbf{Method B}. The proof of low Mach compliance is analogous to that of Method A. For more details, see \cite{barsukow20cgk}.

\begin{remark}
 Observe the presence of denominators of the form $1 + \mathrm{const} \cdot \del_x u$ in both cases, with the constant being (comparable to) the time step. While this is not surprising for method \cite{chalons13}, as it is derived as a Lagrange-Projection method, the presence of the denominator is slightly more intriguing for method \cite{chalons10}, a relaxation solver that does not split the equations and does not treat the compression-sub-system separately. 
\end{remark}

\section{Sequentially-explicit time integration} \label{sec:yee}

Consider the following ``off-diagonal'' system (with $f,g \colon \mathbb R \to \mathbb R \text{ given}$)
\begin{align}
 \del_t \vecc{a}{b} &= \vecc{f(b)}{g(a)}  & a,b &\colon \mathbb R^+_0 \to \mathbb R. \label{eq:examplesystem}
\end{align}
The sequential application of the forward Euler method already uses the value $a^{n+1}$ in the update of $b$:
\begin{align}
 \frac{a^{n+1} - a^n}{\Delta t} &= f(b^n) & 
 \frac{b^{n+1} - b^n}{\Delta t} &= g(a^{n+1}).   \label{eq:triangintro2}
\end{align}

This method is not the forward Euler method. It is explicit, and thus also different from semi-implicit methods of e.g. \cite{guerra86} or \cite{dimarco17}, where similar ideas appear but the final systems remain implicit. The above method of time integration could be called ``leap-frog'', were it not for the many other, and different, methods that already carry this name, such as those in \cite{thomas93}. The name used here, and in \cite{barsukow21yee} therefore is ``sequential explicit''.

One can show that methods using this integration in time have many good properties, such as non-dissipativity (the Fourier modes in the von Neumann stability analysis neither grow, nor decay), and there is a relation to symplectic or energy-conserving time integrators in cases when the PDE at hand possesses a Hamiltonian structure (e.g. those in \cite{remaki99} for the Maxwell equations).

Here, the sequential-explicit integrator shall be merely used for its ability to stabilize central differences. The Euler system can be thought of as being composed of acoustics, advection and compression. It is a lucky coincidence that the acoustic sub-system is of the proper off-diagonal form. As it is only a sub-system of the full Euler equations, a relation to Hamiltonian systems etc. is neither relevant nor required. The advective sub-system, which is not off-diagonal, requires a conventional stabilization mechanism (upwinding). It has also been found in \cite{barsukow21yee} that adding diffusion by way of a denominator as in Lagrange-Projection methods (thus accounting specially for the compressive terms) improves the stability of the overall method. The complete numerical flux (in $x$-direction, say) therefore is, schematically,
\begin{align*}
 f^x_{i+\frac12,j} = \frac{\{\{  \boxed{ \text{flux average}}_{i+\frac12} \}\}_{j\pm\frac12} - \boxed{\text{upwind difference}}_{i+\frac12,j}}{1  + \Delta t \boxed{\text{divergence}}_{i+\frac12,j}}. \label{eq:seqexplicitschematical}
\end{align*}

The upwind difference is that of the advective sub-system only.
In order to solve the acoustic sub-system with sequential-explicit time integration, and the advective one with forward Euler, it is sufficient to update the momentum equations first, and to use the new value of the momentum in the updates of $\rho$ and $e$. This way, the method has the same computational cost as the usual forward Euler method. It has been called Method G in \cite{barsukow21yee}.

It is important to make sure that the cell-based divergence, that arises as $\mathcal O(1)$ term in the energy equation is an average of the edge-based divergence appearing in the denominator, which is why the flux average is a multi-dimensional operator involving averaging in the perpendicular direction. The edge-based divergence, itself, is taken as the average of two node-based divergences, that have already appeared in Section \ref{sec:cgk}. Under this condition one can then prove low Mach compliance (see \cite{barsukow21yee}). This method shall be referred to here as \textbf{Method C}.

\section{Numerical results} \label{sec:numerical}

We aim at comparing the performance of these three very different methods, common properties being merely the explicit nature and the fact that they are conservative methods on Cartesian grids. While method B derives from a Riemann solver, the other two use Riemann solvers at most for sub-systems, with method C essentially being Riemann-solver free. 

\subsection{Kelvin-Helmholtz instability}

The numerical setup chosen for the comparison is a particular choice of a Kelvin-Helmholtz instability from \cite{leidi22,leidi23}. While the proposed methods might have individually been tested already on such setups, here they shall be compared to each other on the same setup for the first time, which shall allow comparison to other low Mach number methods in future.

The setup is studied on a domain $[0,2] \times [-\frac12, \frac12]$ with
\begin{align}
    H(y) &:= \begin{cases} - \sin(\pi \frac{y + \frac14}{w}) & -\frac14 - \frac{w}{2} \leq y < -\frac14 +\frac{w}{2}, \\
     -1 & -\frac14 +\frac{w}{2} \leq y < \frac14 - \frac{w}{2}, \\
     \sin\left(\pi \frac{y - \frac14}{w}\right) & \phantom{-}\frac14 - \frac{w}{2} \leq y < \frac14 + \frac{w}{2}, \\ 
	 1 & \text{else,} \end{cases} 
\end{align}

\begin{align}
	\rho &= \gamma + r H(y), & u &= \mathcal M H(y),  &
	 p &= 1, & v &= \delta \mathcal M \sin(2\pi x).	 
\end{align}

\begin{figure}
 \centering
\includegraphics[width=0.7\textwidth]{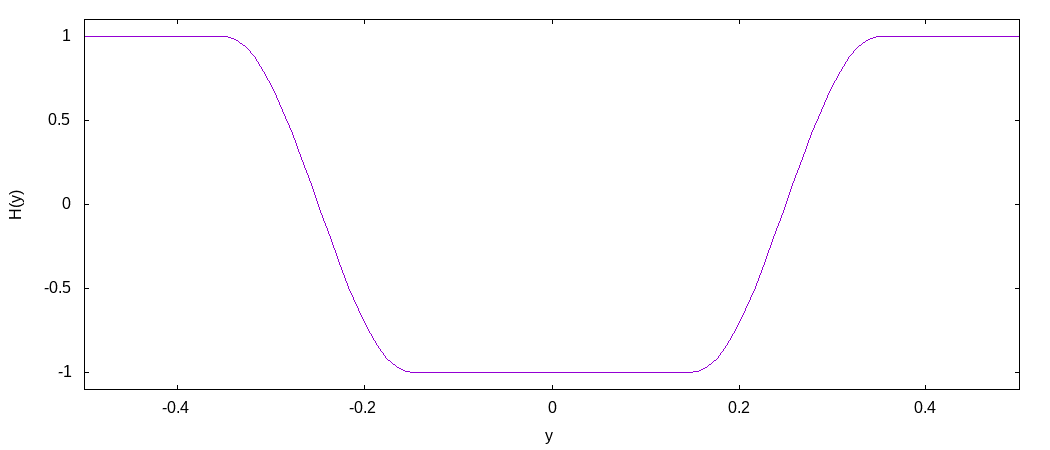}
 \caption{The function $H(y)$, i.e. the vertical flow profile.}
 \label{fig:hfunction}
\end{figure}

In \cite{leidi22}, $r = 0$. Here, we use $r = 10^{-3}$: In incompressible flow, the density acts as a passive scalar, and this way it is possible to visualize the interaction of the layers. Observe that $H \in C^1$. $\mathcal M$ is a free parameter that allows to modify the Mach number of the flow, and this requires a simulation time of $0.8 / \mathcal M$. We choose $\delta = 0.1$ and the width of the shear layer as $w = \frac1{16}$. Periodic boundaries and a CFL number of 0.9 are used, and $\gamma = 1.4$.

\begin{figure}
 \centering
 \includegraphics[width=0.32\textwidth]{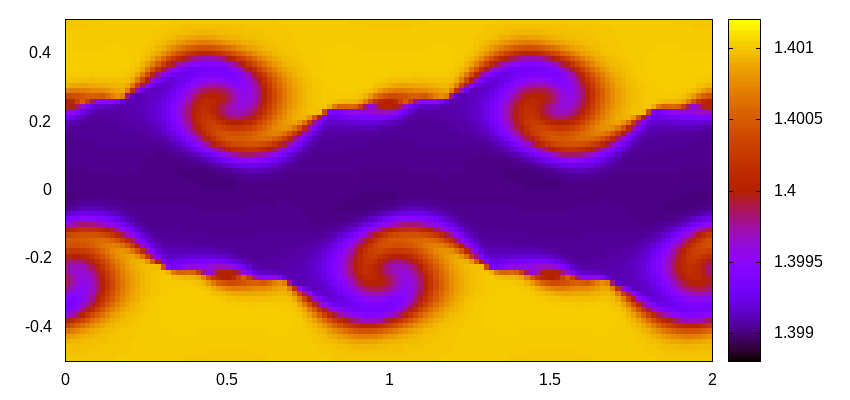}\hfill \includegraphics[width=0.32\textwidth]{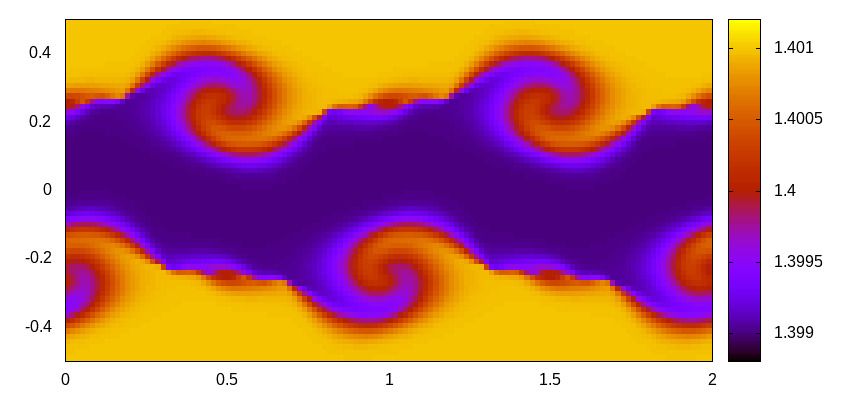}\hfill \includegraphics[width=0.32\textwidth]{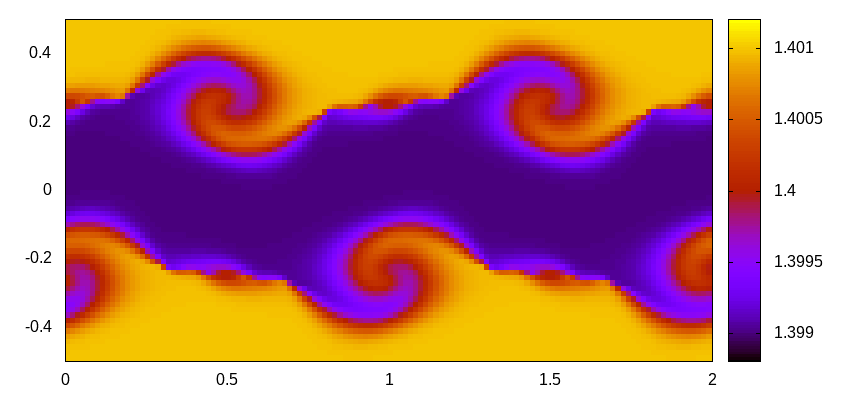} \\
 \includegraphics[width=0.32\textwidth]{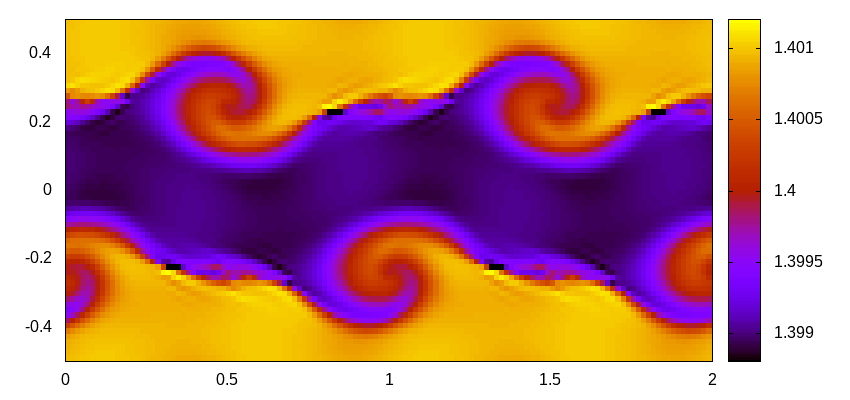}\hfill \includegraphics[width=0.32\textwidth]{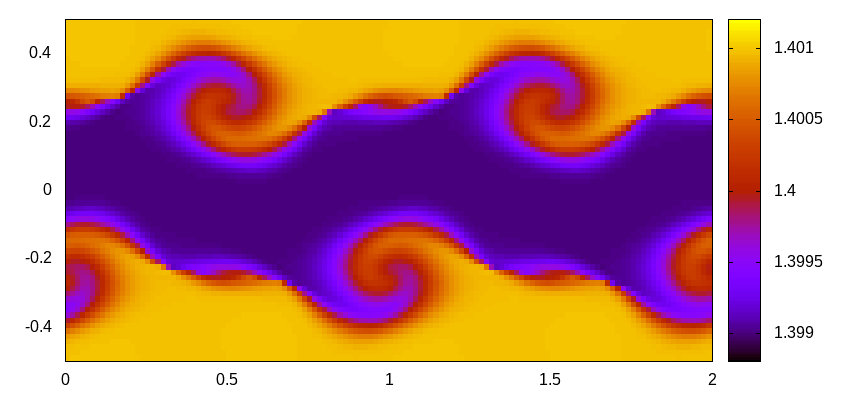}\hfill \includegraphics[width=0.32\textwidth]{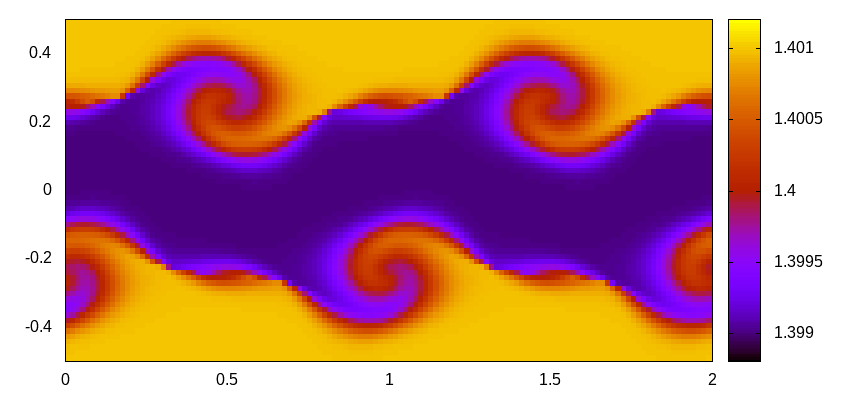} \\
 \includegraphics[width=0.32\textwidth]{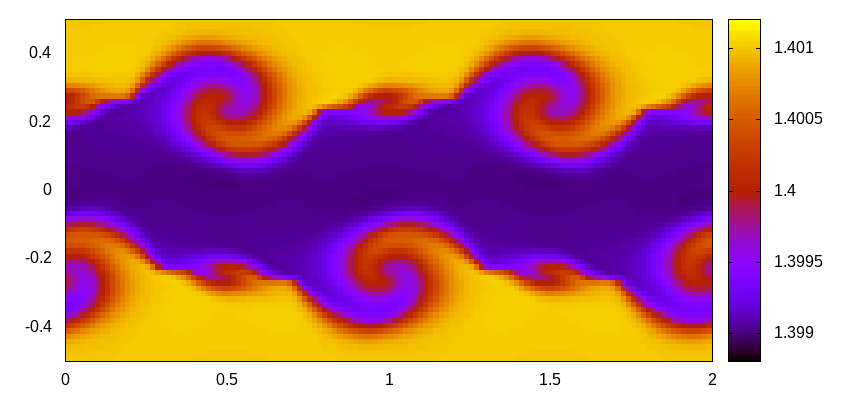}\hfill \includegraphics[width=0.32\textwidth]{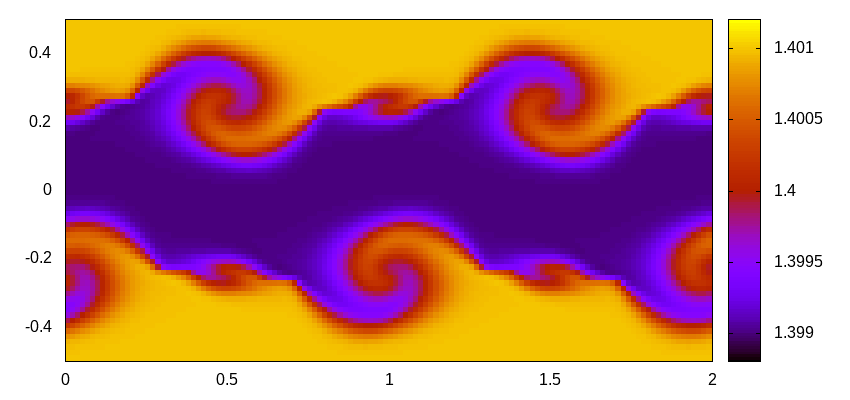}\hfill \includegraphics[width=0.32\textwidth]{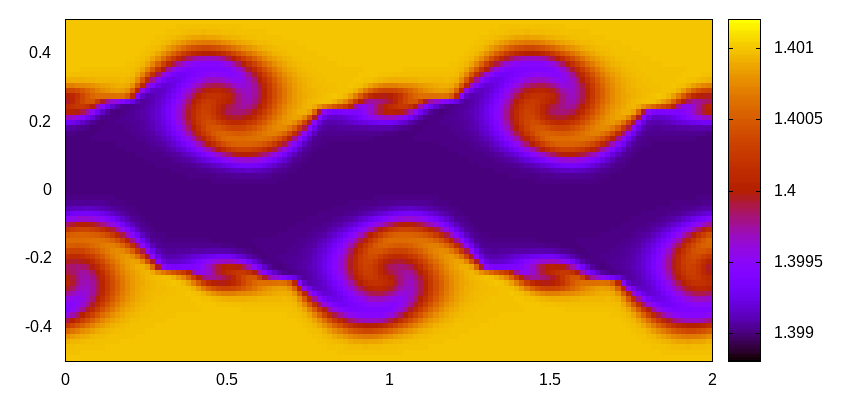} \\
 \caption{Density is color-coded. Grid $128 \times 64$ (all figures using the same setup). Grid $128 \times 64$, $\mathcal M = 10^{-2}, 10^{-3}, 10^{-4}$ (\emph{left to right}). \emph{Top to bottom}: Method A, B, C.}
 \label{fig:tests128vy1}
\end{figure}

\begin{figure}
 \centering
 \includegraphics[width=0.8\textwidth]{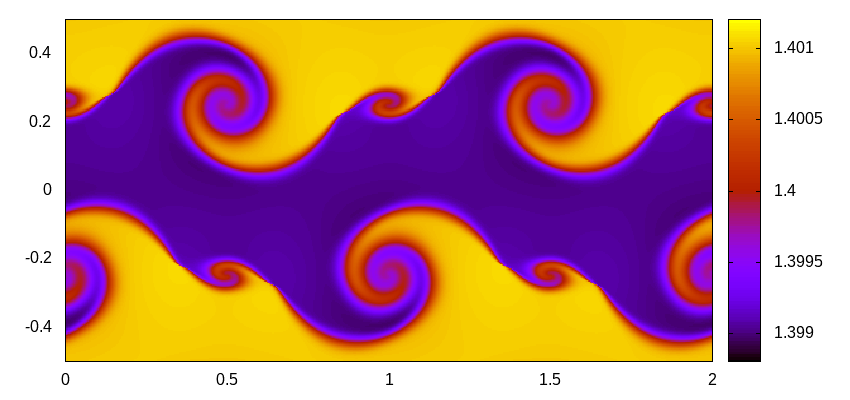} \\
 \includegraphics[width=0.8\textwidth]{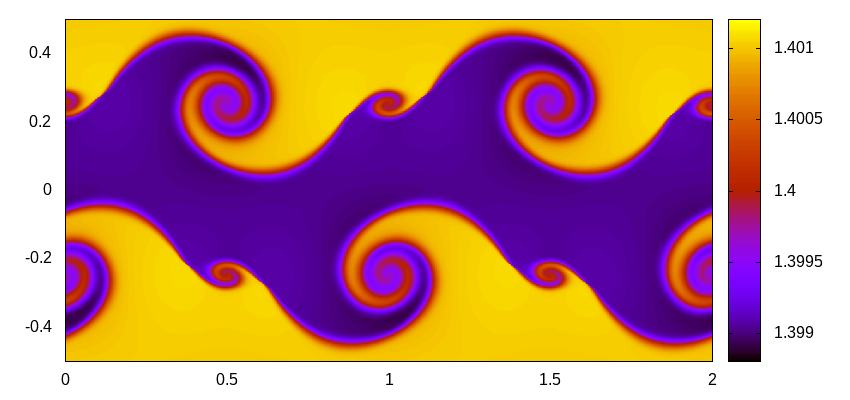} \\
 \includegraphics[width=0.8\textwidth]{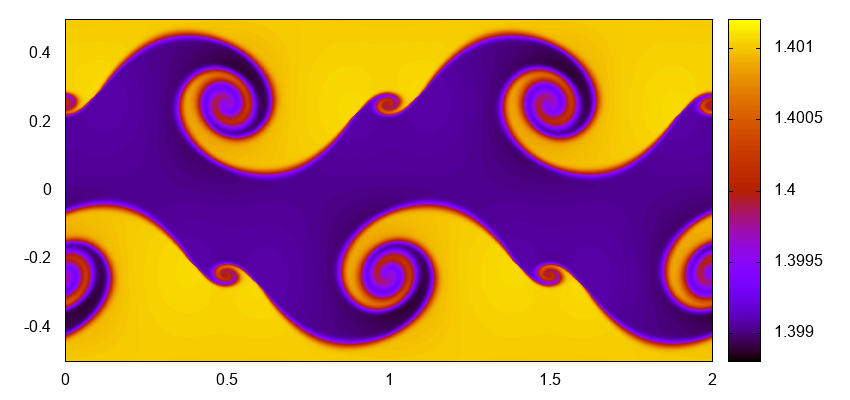} 
 \caption{Convergence study for Method A. Density is color-coded, $\mathcal M = 10^{-2}$. \emph{Top to bottom}: Grid $512 \times 256$, $1024 \times 512$, $2048 \times 1024$.}
 \label{fig:testsm2convergence}
\end{figure}

\begin{figure}
 \centering
\includegraphics[width=0.7\textwidth]{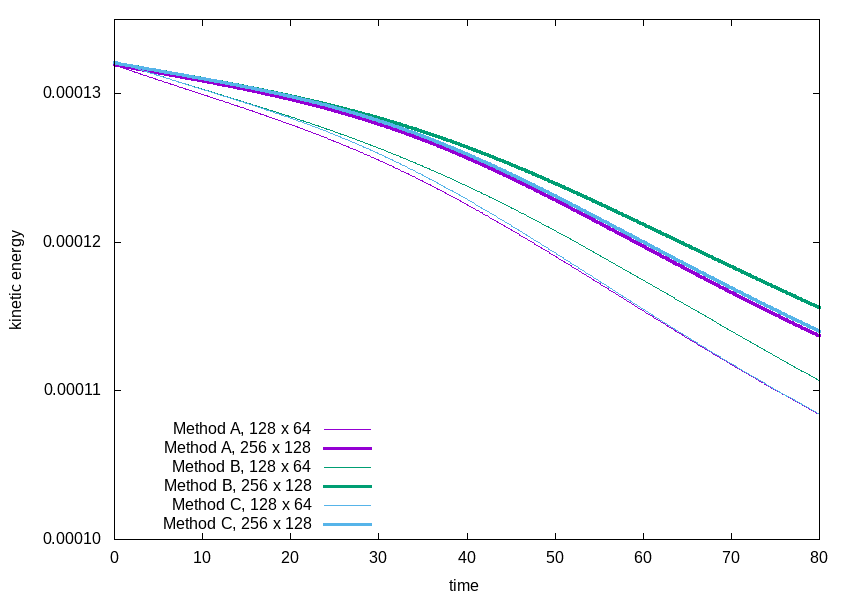}
 \caption{The evolution of kinetic energy for the three methods ($\mathcal M = 10^{-2}$). In the incompressible limit, kinetic energy is conserved; the observed decay therefore is entirely due to the numerical diffusion. As the curves do not change for $\mathcal M = 10^{-3}$ (up to a rescaling of the energy by $\mathcal M^{-2}$ and of time by $\mathcal M$) (not shown), one can even conclude that the decay is due to the numerical diffusion of the advective operator.}
 \label{fig:kinen}
\end{figure}

One observes (Figure \ref{fig:tests128vy1}) that indeed the results are independent of the Mach number, which is the essence of low Mach compliance. Despite the very different strategies of derivation, the methods yield very similar results. One observes the appearance of artefacts (most pronounced for Method C), where a part of the shear layer seems not to become unstable in the proper way. They develop into vortices upon further grid refinement (Figure \ref{fig:testsm2convergence}), but the vortices keep decreasing in size and are purely numerical. They tend not to appear upon usage of a method of second or higher order of accuracy (\cite{leidi23}). On a grid of $128 \times 64$ ($256 \times 128$) one finds (see Figure \ref{fig:kinen}) the total kinetic energy to have decayed by the following amounts by the end of the simulations: Method A 17.8\% (13.8\%), Method B 16.3\% (12.5\%), Method C 18.0\% (13.7\%). These values are indistinguishable for $\mathcal M=10^{-2}$ and $\mathcal M = 10^{-3}$, which means that the decay is due to the advective diffusion. Methods A and C are designed to contain more of it, which might explain the slightly higher loss of kinetic energy.

\subsection{Radial Riemann problem}

\begin{figure}
 \centering
 \includegraphics[width=0.49\textwidth]{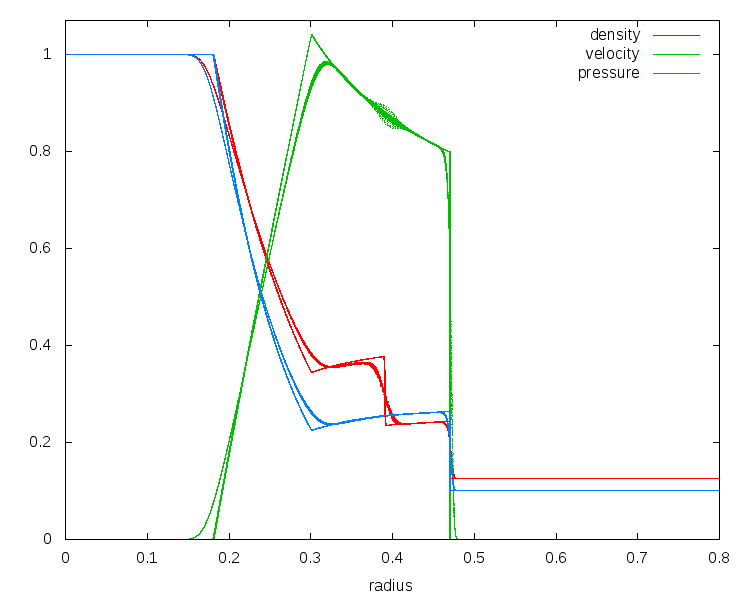} \includegraphics[width=0.49\textwidth]{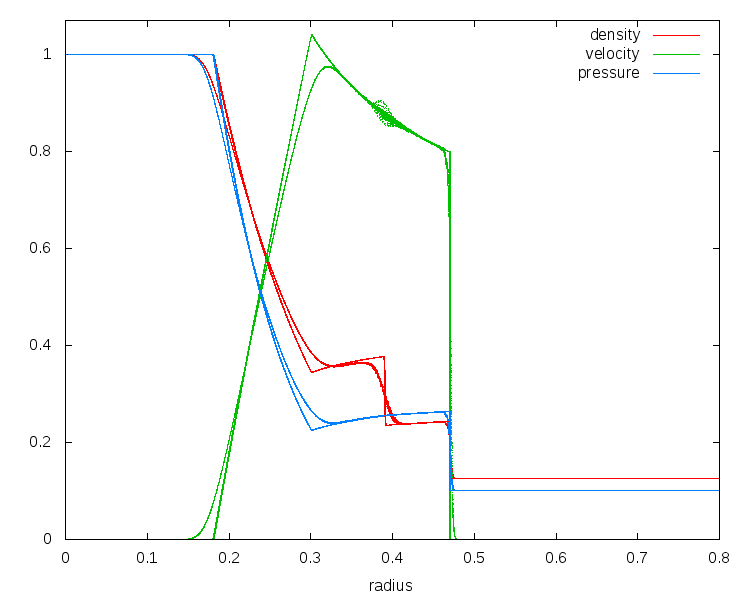}  \includegraphics[width=0.49\textwidth]{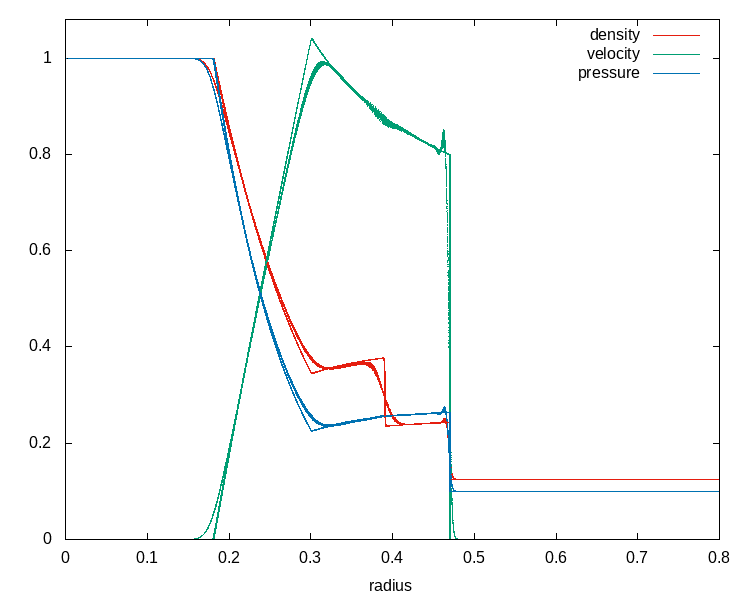} 
 \caption{Scatter plots of the numerical solution at $t=0.1$ of a radial Riemann problem. \emph{Top row}: Method A (\emph{left}) and Method B (\emph{right}) (Figures from \cite{barsukow20cgk}.) \emph{Bottom row}: Method C. The solid lines show the solution obtained by solving just the radial, one-dimensional equations on a very fine grid.}
 \label{fig:sodshock}
\end{figure}

The all-speed property of the three methods shall finally be assessed using a radial version of the Sod shock tube. Methods A and B are multi-dimensional extensions of stable one-dimensional methods, and this test allows to see how they cope with a multi-dimensional setting. This setup is also interesting for Method C, as it is Riemann-solver-free. The initial data are $\rho = 0.125$, $p = 0.1$ outside of a disc with radius 0.3, and $\rho = 1$, $p = 1$ inside. The velocity is initially zero everywhere. Figure \ref{fig:sodshock} shows radial scatter plots for the three methods on a $500\times 500$ grid at time $t=0.1$. One observes that Methods A and B give very similar results, with a slight scatter (i.e. a defect of rotational symmetry) around the contant wave. Method C shows a small oscillation around the shock.

\section{Conclusions}

This paper presents numerical results of three recently developed truly-multidi\-men\-sional numerical methods for the Euler equations. They achieve the all-speed property without ad hoc fixes, by including the divergence operator in the numerical diffusion instead of derivatives of individual velocity components. Despite being first-order accurate they resolve details of a complex setup of a Kelvin-Helmholtz instability even on coarse grids and show asymptotically no Mach number dependence. They, however, display artefacts in the form of secondary vortices that are purely numerical. Future research will focus on the precise origin of the artefacts and on extensions of these methods to higher order, which will also include limiting.

\section{Acknowledgements}

The author thanks Giovanni Leidi (HITS, Heidelberg) for inspiring discussions.

\bibliographystyle{alpha}
\newcommand{\etalchar}[1]{$^{#1}$}

\end{document}